\documentclass[a4paper,10pt]{article}
\usepackage{amsmath,amssymb}
\usepackage[latin1]{inputenc}

\newcommand{\prelie}{\operatorname{PreLie}}

\newcommand{\lie}{\operatorname{Lie}}

\newcommand{\PL}{\mathtt{pL}}
\newcommand{\ad}{\mathtt{ad}}
\newcommand{\Ks}{\mathtt{K}}
\newcommand{\kk}{\mathbb{K}}
\newcommand{\sym}{\mathfrak{S}}
\newcommand{\uW}{\mathsf{U}(W)}
\newcommand{\sous}{\curvearrowleft}
\newcommand{\Free}{\operatorname{Free}}

\newtheorem{theorem}{Theorem}[section] 
\newtheorem{proposition}[theorem]{Proposition} 
 
\newtheorem{corollary}[theorem]{Corollary} 
\newtheorem{lemma}[theorem]{Lemma}

\newenvironment{proof}{\begin{trivlist}\item{\bf{Proof.}}}
  {\hfill\rule{2mm}{2mm}\end{trivlist}}

\title{Free pre-Lie algebras are free as Lie algebras} 

\author{F. Chapoton}
\date{\today}

\begin{document}
\maketitle

\begin{abstract}
  We prove that free pre-Lie algebras, when considered as Lie
  algebras, are free. Working in the category of $\sym$-modules, we
  define a natural filtration on the space of generators. We also
  relate the symmetric group action on generators with the structure
  of the anticyclic PreLie operad.
\end{abstract}


\section*{Introduction}

A pre-Lie algebra is a vector space $V$ endowed with a bilinear map
$\sous\, : V\otimes V \to V$ such that
\begin{equation}
  \label{axiomePL}
  (x \sous y) \sous z -  x \sous (y \sous z)=  (x \sous z) \sous y -  x \sous (z \sous y),
\end{equation}
for all $x,y,z\in V$. This kind of algebra has been used for a long
time in various areas, see \cite{burde} for a survey. In geometry,
there is a pre-Lie product on the space of vector fields of any
variety endowed with an affine structure. In algebra, this kind of product
appears for instance in relation with deformation theory, operads and
vertex algebras.

In a previous article \cite{prelie}, the free pre-Lie algebras have been
completely described in terms of rooted trees. The language of operads
was a convenient setting for this.

In any pre-Lie algebra, the bracket  $[x,y]=x\sous y-y \sous
x$ satisfies the Jacobi identity and therefore defines a Lie algebra on
the same vector space.

The main goal of this article is to prove that free pre-Lie algebras
are free Lie algebras for the Lie bracket $[\,,\,]$. This result has
been obtained before with different methods by Foissy \cite{foissy}.
We give here a new proof using a spectral sequence and obtain a
natural filtration on the space of generators. 

\medskip

The article is written using the language of $\sym$-modules. We will
briefly recall the main features of this theory. This is the natural
category for operads and provides a clean way to deal with free
algebras. 

The main result is the fact that the $\sym$-module for pre-Lie
algebras is isomorphic (as a Lie algebra in the category of
$\sym$-modules) to the composition of the $\sym$-module for Lie
algebras with a $\sym$-module $X$. This implies the statement above
about free pre-Lie algebras. We also give a description of the
$\sym$-module $X$ and a relation between $X$ and the anticyclic
structure on the PreLie operad.

Thanks to Jérome Germoni and Muriel Livernet for their interest in
this work and useful suggestions.

\section{Notations and $\sym$-modules}

Let $\kk$ be a fixed ground field of characteristic zero.

We will work in the Abelian category of $\sym$-modules over $\kk$ or
sometimes with complexes of $\sym$-modules over $\kk$. Recall that a
$\sym$-module $P$ can either be seen as a functor $P$ from the
category of finite sets and bijections to the category of vector
spaces over $\kk$ or as a sequence $(P(n))_{n \geq 0}$ of
$\kk$-modules over the symmetric groups $\sym_n$. We will use freely
one or the other of these equivalent definitions.

The category of $\sym$-modules is symmetric monoidal for the following
tensor product:
\begin{equation}
  (F\otimes G)(I)=\sqcup_{I=J \sqcup K} F(J) \times G(K),
\end{equation}
where $I$ is a finite set.

There is another monoidal structure on $\sym$-modules which is
non-symmetric and defined by
\begin{equation}
(F \circ G)(I) =\sqcup_{\simeq} F(I/{\simeq}) \times \prod_{J
\in I/\simeq} G(J),
\end{equation}
where $I$ is a finite set and $\simeq$ runs over the set of equivalence relations on $I$.

An operad $Q$ is essentially a monoid in the monoidal category of
$\sym$-modules with $\circ$ as tensor product, \textit{i.e.} the data of a map
from $Q \circ Q$ to $Q$ which is associative.

We will consider $\sym$-modules endowed with various kinds of
algebraic structures. One can define associative algebras, commutative
algebras, exterior algebras, Lie algebras, Hopf algebras and pre-Lie
algebras in the category of $\sym$-modules by using the symmetric
monoidal product $\otimes$ and the usual diagrammatic definitions.
This works in particular for algebras over any operad $Q$. This kind
of object is sometimes called a twisted algebra or a left-module over
the corresponding operad. We will simply call them algebras.

If $Q$ is an operad and $P$ is an $\sym$-module, a structure of
$Q$-algebra on $P$ is given by a map from $Q \circ P$ to $P$.

If $Q$ is an operad and $P$ is an $\sym$-module, the $\sym$-module $Q
\circ P$ is the free $Q$-algebra on $P$. The map from $Q \circ Q \circ
P$ to $Q\circ P$ is deduced from the map from $Q \circ Q$ to $Q$
defining the operad $Q$.

\medskip

For instance, let us consider the following $\sym$-modules.

Let $\lambda$ be a partition of an integer $n$. Then $\lambda$ is
associated with an irreducible representation of the symmetric group
$\sym_n$, that will be denoted $S^{\lambda}$. This can be seen as a
$\sym$-module concentrated in degree $n$. As a special case, $S^1$ is
the trivial module over $\sym_1$. Let $S$ be the direct sum of all
$\sym$-modules $S^n$ corresponding to trivial modules over the
symmetric groups.

Let $T$ be the $\sym$-module defined as the direct sum of all the regular
representations $T^n$ over the symmetric groups.

Let $\Lambda$ be the direct sum of all $\sym$-modules $\Lambda^n$
corresponding to alternating modules over the symmetric groups.

There are morphisms of $\sym$-modules $S\otimes S\to S$, $T\otimes
T\to T$ and $\Lambda \otimes \Lambda \to \Lambda$, which define
associative algebras in the category of $\sym$-modules. These algebras
are essentially versions of the polynomial algebra, the tensor algebra
and the exterior algebra.

For example, if $I$ is a finite set, $T(I)$ is the vector space
spanned by total orders on $I$. One can concatenate total orders on
two sets $I$ and $J$ to get a total order on $I \sqcup J$. This
defines the map from $T \otimes T$ to $T$.

\medskip

Let $P$ be an  $\sym$-module. The generating series $f_P$ associated
with $P$ is defined by
\begin{equation}
  f_P=\sum_{n\geq 0} \dim(P(\{1,\dots,n\})) \frac{x^n}{n!}.
\end{equation}
There is a more refined object, which is a symmetric function $Z_P$
recording the action of the symmetric groups:
\begin{equation}
  Z_P=\sum_{n\geq 0} \chi(P(\{1,\dots,n\})),
\end{equation}
where $\chi$ is the image of the character of the $\sym_n$-module
$P(\{1,\dots,n\})$ in the ring of symmetric functions. This symmetric
function is a complete invariant, in the sense that two $\sym$-modules
are isomorphic if and only if they share the same symmetric function.



 

For more information on $\sym$-modules and operads, see for instance
\cite{muriel,renaissance}.

\subsection{Relation with usual algebra}

Each $\sym$-module $P$
defines a functor from vector spaces to vector spaces as follows:
\begin{equation}
  P(V)=\oplus_{n \geq 0} P(n)\otimes_{\sym_n} V^{\otimes n}.
\end{equation}

The tensor product $\otimes$ of $\sym$-modules has the following
property:
\begin{equation}
  (F\otimes G)(V)\simeq F(V)\otimes G(V),
\end{equation}
where $V$ is a vector space.

The tensor product $\circ$ of $\sym$-modules has the following
property:
\begin{equation}
  (F\circ G)(V)\simeq F(G(V)),
\end{equation}
where $V$ is a vector space.

If the $\sym$-module $P$ is an operad and $V$ is a vector space, then
$P(V)$ is the free $P$-algebra on $V$.

Then, if $P$ is a $\sym$-module with such an algebraic structure,  \textit{i.e.}
a $Q$-algebra for some operad $Q$, the vector space $P(V)$ is endowed
with the corresponding structure of $Q$-algebra in the category of
vector spaces.

For instance, the classical algebra structures on $S(V)$, $T(V)$ and
$\Lambda(V)$ (respectively the polynomial algebra, the tensor algebra
and the exterior algebra) come from the morphisms of $\sym$-modules
$S\otimes S\to S$, $T\otimes T\to T$ and $\Lambda \otimes \Lambda \to
\Lambda$ that were introduced above.

\section{The pre-Lie $\sym$-module $W$}

A rooted tree on a set $I$ is a connected and simply-connected graph
with vertex set $I$ together with a distinguished element of $I$
called the root. Note that a rooted tree can be canonically decomposed
into its roots and a set of rooted trees obtained when the root is
removed. Let $W$ be the $\sym$-module which maps a finite set $I$ to
the vector space spanned by the set of rooted trees on $I$. 

It has been shown in \cite{prelie} that $W$ is the operad describing
pre-Lie algebras and in particular that $W(V)$ is the free pre-Lie
algebra on a vector space $V$. There is a natural morphism $\sous : W
\otimes W \to W$ corresponding to the PreLie product. This means that
$W$ is a pre-Lie algebra in the category of $\sym$-modules.

Let us recall briefly the definition of the map $\sous$.

Given a rooted tree $S$ on a finite set $I$ and a rooted tree $T$ on a
finite set $J$, one can define $S \sous T$ (which is a sum of rooted
trees on the set $I \sqcup J$) as follows: consider the disjoint union
of $S$ and $T$, then add an edge between the root of $T$ and a vertex
of $S$ and sum with respect to the chosen vertex of $S$. The root is
taken to be that of $S$.

As there exists a morphism from the $\lie$ operad to the $\prelie$
operad given by 
\begin{equation}
  [x,y]=x\sous y-y\sous x,
\end{equation}
there is a morphism $[\,,\,] : W \otimes W \to W$ which makes $W$ into
a Lie algebra in the category of $\sym$-modules.

The universal enveloping algebra of a Lie algebra is a well-defined
associative algebra in the category of $\sym$-modules, which has most
of the classical properties of the usual construction for vector
spaces. Let $\mathsf{U}(W)$ be the universal enveloping algebra of
$W$.

In the sequel, we will call \textbf{$W$-module} a module over the
Lie algebra $W$ or equivalently a right module over
$\mathsf{U}(W)$. This notion is defined in the obvious way in the
symmetric monoidal category of $\sym$-modules.

There are two distinct $W$-module structures on the $\sym$-module $W$.
The first one is given by the adjoint action of the Lie algebra $W$ on
itself. It will be denoted by $W_{\ad}$. The other one is given by the
PreLie product $\sous$ and will be denoted by $W_{\PL}$. The fact that
the map $\sous$ is a right action results from the PreLie axiom
(\ref{axiomePL}).

\bigskip

Let us now recall results from \cite[Th. 3.3]{prelie}. There is an
isomorphism $\psi$ of $W$-modules between $W_{\PL}$ and the free
$\uW$-module on the $\sym$-module $S^1$. By the description of $W$ in
terms of rooted trees and the decomposition of a rooted tree into its
root and its set of subtrees, the module $W_{\PL}$ is isomorphic as an
$\sym$-module to $S^1 \otimes S \circ W$. The $\sym$-module $S \circ
W$ is spanned by forests of rooted trees. The isomorphism $\psi$ can
be written as $\operatorname{Id} \otimes \phi$ between $S^1 \otimes S
\circ W$ and $S^1 \otimes \uW$. This defines an isomorphism $\phi$ of
$W$-modules from $S \circ W$ to $\uW$, where the action
$\triangleleft$ of $W$ on $S \circ W$ is deduced from the product
$\sous$ and is given by
\begin{equation}
  \label{action_totale}
  (t_1 t_2 \dots t_k )\triangleleft t =  t_1 t_2 \dots t_k
  t+\sum_{i=1}^{k} t_1 t_2 \dots t_{i-1} (t_i \sous t) t_{i+1} \dots t_k,
\end{equation}
where the $t$'s stand for some rooted trees.

We will from now on identify by the mean of $\phi$ the $W$-module
$\uW$ with the $W$-module $S \circ W$ with this action
$\triangleleft$. One can see from the explicit shape of the action
$\triangleleft$ that the $W$-module $\uW$ has a decreasing filtration
by the number of connected components of the forest.

The associated graded module is given by the action
\begin{equation}
  \label{actionPL}
  (t_1 t_2 \dots t_k ) \sous t = \sum_{i=1}^{k} t_1 t_2 \dots t_{i-1} (t_i \sous t) t_{i+1} \dots t_k,
\end{equation}
where we have slightly abused notation by using the symbol $\sous$ for
the action. This is easily seen to be the natural $W$-module structure
on $S \circ W_{\PL}$ obtained by extending the $W$-module $W_{\PL}$ by
derivation. This is also the symmetric algebra on the $W$-module
$W_{\PL}$.

One can see that there is only one other term in the filtered action
(\ref{action_totale}), which is of degree $1$ with respect to the
graduation by the number of connected components and is just a
product:
\begin{equation}
  \label{actionK}
  t_1 t_2 \dots t_k \otimes t \mapsto  t_1 t_2 \dots t_k  t.
\end{equation}

\section{Two spectral sequences}

Let us consider the usual reduced complex computing the homology of
the Lie algebra $W$ with coefficients in the $W$-module $\uW$. This is
the tensor product of the exterior algebra on $W$ with the module
$\uW \simeq S \circ W$.

As an $\sym$-module, this complex is $(S \circ W) \otimes (\Lambda
\otimes W)$. The differential $\partial$ is the usual
Chevalley-Eilenberg map, which uses both the bracket or the action as a
contraction:
\begin{multline}
  \partial(x_1 x_2 \dots x_k \otimes y_1 \wedge y_2 \wedge \dots\wedge
  y_\ell)
= \sum_{j=1}^{\ell} \pm (x_1 x_2 \dots x_k) \triangleleft y_j \otimes  y_1 \wedge \dots \wedge \widehat{y_j}
\wedge \dots \wedge y_\ell  + 
\\ \sum_{1\leq i<j\leq \ell} \pm x_1 x_2 \dots x_k
\otimes [y_i,y_j] \wedge y_1 \wedge \dots \wedge \widehat{y_{i}}\wedge \dots \wedge 
\widehat{y_{j}} \wedge \dots 
\wedge y_\ell,
\end{multline}
where the signs are given by the Koszul sign rule.

But $\uW$ is a free $W$-module by definition, hence the homology is
concentrated in homological degree $0$ and is given by the
$\sym$-module $S^0$.

Let us now use the filtration on $\uW$ to define two spectral
sequences computing the same homology. In fact, we will first define a
bicomplex and then consider its two associated spectral sequences.

We have to introduce a triple grading on the complex $(S \circ W)
\otimes (\Lambda \otimes W)$. Let us denote by $n$ the internal degree
of the $\sym$-module (seen as a collection of modules over $\sym_n$),
by $p$ the degree with respect to the graduation of the symmetric
algebra and by $q$ the homological degree with respect to the
graduation in the exterior algebra. As $W$ has no component with
$n=0$, one has $p\geq 0$, $q\geq 0$ et $p+q \leq n$.

Let $r$ be $n-p-q$. We will use the triple grading by $(n,p,r)$. The
differential $\partial$ is of degree $0$ with respect to the first
grading. Hence one can consider separately each part of fixed first
degree $n$.

The differential $\partial$ on  $(S \circ W)
\otimes (\Lambda \otimes W)$ decomposes into two pieces according to
the decomposition of the action $\triangleleft$ into the action
$\sous$ coming from (\ref{actionPL}) plus
another term coming from (\ref{actionK}).

The first part is defined as follows:
\begin{multline}
  \partial_{\PL} (x_1 x_2 \dots x_k \otimes y_1 \wedge y_2 \wedge \dots\wedge
  y_\ell)
= \sum_{j=1}^{\ell} \pm (x_1 x_2 \dots x_k) \sous y_j \otimes  y_1 \wedge \dots \wedge \widehat{y_j}
\wedge \dots \wedge y_\ell   + \\
 \sum_{1\leq i<j\leq \ell} \pm x_1 x_2 \dots x_k
\otimes [y_i,y_j] \wedge y_1 \wedge \dots \wedge \widehat{y_{i}}\wedge \dots \wedge 
\widehat{y_{j}} \wedge \dots 
\wedge y_\ell.
\end{multline}
One can recognize the differential $\partial_{\PL}$ of degree
$(0,0,1)$ computing the homology of the Lie algebra $W$ with
coefficients in the graded $W$-module $(S \circ W_{\PL},\sous)$.

The remaining terms of the differential are
\begin{equation}
  \partial_{\Ks} (x_1 x_2 \dots x_k \otimes y_1 \wedge y_2 \wedge \dots\wedge
  y_\ell) 
= \sum_{j=1}^{\ell} \pm x_1 x_2 \dots x_k y_j \otimes y_1 \wedge \dots \wedge \widehat{y_j}
\wedge \dots \wedge y_\ell. 
\end{equation}
The map $\partial_{\Ks}$ has degree $(0,1,0)$ and does no longer use
the bracket of the Lie algebra $W$. This is nothing but the
differential in the Koszul complex relating the exterior algebra
$\Lambda \circ W$ on $W$ and the symmetric algebra $S \circ W$ on $W$
\cite{priddy}.

The two differentials $\partial_{\PL}$ and $\partial_{\Ks}$ are of
degree $(0,0,1)$ and $(0,1,0)$ and their sum is also a differential.
Hence they define a bicomplex and one can consider the two spectral
sequences associated with this bicomplex.

\begin{proposition}
  The spectral sequence beginning with $\partial_{\Ks}$ degenerates at
  first step.
\end{proposition}

\begin{proof}
  As it is known that the exterior and symmetric algebras are Koszul
  dual of each other and Koszul, the homology of $\partial_{\Ks}$ is
  concentrated in homological degree $0$ and is given by $S^0$.
\end{proof}

Before studying the other spectral sequence, one needs the two
following results. Let us first recall a classical lemma:
\begin{lemma}
  \label{jantzen}
  Let $A$ be a Hopf algebra and let $N$ be a right $A$-module. Then
  $N\otimes A$ is isomorphic as a right $A$-module to the free right $A$-module
  generated by $N$.
\end{lemma}
\begin{proof}
  The argument uses the antipode of $A$ to define an isomorphism, see for
  instance \cite[\S 3.6 and \S 3.7]{jantzen}.
\end{proof}

We will also need the following property of some $W$-modules.

\begin{proposition}
  \label{projective}
  Let $\lambda$ be a non-empty partition. Then $S^\lambda \circ W_{\PL}$ is
  a projective $\uW$-module.
\end{proposition}

\begin{proof}
  Recall that $W_{\PL}$ is isomorphic to the free $\uW$-module $S^1
  \otimes \uW$.

  Hence, for each integer $k\geq 1$, the $W$-module $T^k \circ
  W_{\PL}$ is isomorphic to the module $T^k \otimes (T^k \circ \uW)$,
  where $W$ acts on the right factor only. Here, we have used the
  property of the $\sym$-module $T$ that $T \circ (A \otimes B) \simeq
  (T \circ A) \otimes (T \circ B)$.

  As $\uW$ is a Hopf algebra, one can apply Lemma \ref{jantzen}, which
  implies that $T^k \circ \uW$ is a free $\uW$-module. It follows that
   $T^k \circ W_{\PL}$ is also a free $\uW$-module.

   Now, for each partition $\lambda$, $S^\lambda$ is usually defined
   as a direct factor of $T^{|\lambda|}$, where $|\lambda| $ is the
   size of $\lambda$. So $S^\lambda \circ W$ is a direct factor of
   $T^{|\lambda|} \circ W$. As $T^{|\lambda|} \circ W$ is a free $\uW$-module,
   $S^\lambda \circ W$ is a projective $\uW$-module.
\end{proof}

Remark: it may be that all the $\uW$-modules $S^\lambda \circ W$ are
in fact free.

\bigskip

Let us now go back to the bicomplex. To illustrate the computation of
the horizontal spectral sequence starting with $\partial_{\PL}$, we
will draw the component of fixed first degree $n$ of the bicomplex in
the first quadrant, with the $p$ grading increasing from bottom to top
and the $r$ grading from left to right.

We have the following description of the homology with respect to
$\partial_{\PL}$.

The bottom line ($p=0$) of the bicomplex is the complex $S^0 \otimes
(\Lambda \circ W)$ computing the homology of the Lie algebra $W$ with
coefficients in the trivial module. This is what we would like to
compute.

The other lines ($p>0$) of the bicomplex are the complexes $(S^p \circ
W_{\PL}) \otimes (\Lambda \circ W)$ computing the homology of the Lie
algebra $W$ with coefficients in the modules $(S^p \circ W_{\PL})$. As
we know that these modules are projective by Prop. \ref{projective},
the homology is concentrated in degree $q=0$  \textit{i.e.} $r=n-p$.

So the first step of the spectral sequence looks as follows:
\begin{equation}
  \begin{matrix}
    * & & & & & \\
    0 & * & & & & \\
    0 & 0 & *& & & \\
    0 & 0 & 0& *& & \\
    0 & 0 & 0 & 0 & *& \\
    * & * & * & * & * & *,\\
 \end{matrix}
\end{equation}
where a $*$ means a possibly non-zero homology group.

As the spectral sequence converges to $S^0$ because the homology of
the total complex is $S^0$, one deduces from the shape above that the
homology of the bottom line is concentrated in degree $q=1$  \textit{i.e.}
$r=n-1$. So in fact, the first step looks like that:
\begin{equation}
  \label{step1}
  \begin{matrix}
    * & & & & & \\
    0 & * & & & & \\
    0 & 0 & *& & & \\
    0 & 0 & 0& *& & \\
    0 & 0 & 0 & 0 & *& \\
    0 & 0 & 0 & 0 & * & 0.\\
 \end{matrix}
\end{equation}

Let us denote by $H(n,p,r)$ these homology groups. We will denote by
$X$ the $\sym$-module corresponding to the collection of
$\sym_n$-modules $H(n,0,n-1)$ for $n\geq 1$. Note that one does not
include the degree $0$ component in this definition.

\begin{proposition}
  Let $n\geq 1$. The dimension of $H(n,0,n-1)$ is $(n-1)^{n-1}$. The
  dimension of $H(n,p,n-p)$ is $\binom{n}{p}(p-1)(n-1)^{n-p-1}$ for $1
  \leq p \leq n$. There is a filtration on the unique non-vanishing
  homology group $H(n,0,n-1)$ of the bottom line whose graded pieces
  are isomorphic to the homology groups $H(n,p,n-p)$ of the other
  lines.
\end{proposition}

\begin{proof}
  Using once again the fact that the spectral sequence converges to
  $S^0$, one can see that the component of first degree $n$ of the
  horizontal spectral sequence cannot degenerate before the nth step
  and that the successive pages of the spectral sequence provide the
  expected filtration on $H(n,0,n-1)$.

  Let $f_W$ be the generating series associated with $W$:
  \begin{equation}
  \label{taylorW}
    f_W=\sum_{n\geq 1} n^{n-1} \frac{x^n}{n!}.
  \end{equation}

  The function $-f_W(-x)$ is usually called the Lambert W function
  \cite{corless}.
  
  Let us consider now the graded $\sym$-module $S \circ W$. We introduce
  a variable $s$ to encode the $p$ grading. The associated generating
  series is
  \begin{equation}
    f_{S \circ W}=e^{s f_W}.
  \end{equation}
  
  For the graded exterior $\sym$-module $\Lambda \circ W$, we
  introduce a variable $-t$ to encode the $q$ grading. The minus sign
  in front of $t$ is convenient here ; specialization at $t=1$ gives
  the Euler characteristic. The associated generating series is
  \begin{equation}
    f_{\Lambda \circ W}=e^{-t f_W}.
  \end{equation}

  The associated generating series for the bicomplex is then given by
  \begin{equation}
    f_{(S \circ W)
      \otimes (\Lambda \otimes W)} = e^{s f_W}\,e^{-t f_W}=e^{(s-t)f_W}=1+(s-t)\sum_{n\geq 1}(n+s-t)^{n-1}\frac{x^n}{n!}.
  \end{equation}
  Here the Taylor expansion is a classical result, see for example
  \cite[Formula (2.36)]{corless}.
  
  Taking the horizontal Euler characteristic is given by substituting
  $t=1$:
  \begin{equation}
    \label{euler_series}
    e^{s f_W}e^{- f_W}=e^{(s-1)f_W}=1+(s-1)\sum_{n\geq 1}(n+s-1)^{n-1}\frac{x^n}{n!}.
  \end{equation}
  
  As we known by (\ref{step1}) where the horizontal homology of the
  bicomplex is concentrated, computing the Euler characteristic is
  enough to get the dimension of the homology.

  Let us first compute the constant term of (\ref{euler_series}) with
  respect to $s$. One gets
  \begin{equation}
    1-\sum_{n\geq 1}(n-1)^{n-1}\frac{x^n}{n!}.
  \end{equation}
  Hence the dimension of $H(n,0,n-1)$ is $(n-1)^{n-1}$ as expected.

  One can also easily compute the coefficient of $s^p$ for $p>0$ and
  get the expected formula for the dimension of $H(n,p,n-p)$.
\end{proof}

Remark: the filtration on $H(n,0,n-1)$ with quotients $H(n,p,n-p)$ for
$1\leq p\leq n$ gives a nice interpretation of the classical identity
\begin{equation}
  (n-1)^{n-1}=\sum_{p=1}^{n} \binom{n}{p}(p-1)(n-1)^{n-p-1},
\end{equation}
which can be found for instance in \cite[Prop. 2]{chauve} and goes
back to Cayley \cite{cayley}.

\section{Description of $X$}

One can get not just the dimensions of $H(n,0,n-1)$ but a description
of the action of the symmetric groups on $X$.

In \cite[Prop. 7.2]{hypertree}, the symmetric function $Z_{\Lambda
  \circ W}$ was computed, with a parameter $-t$ accounting for the
cohomological grading. Putting $t=1$ in the formula there and removing
the constant term, one finds the symmetric function for $X$.

Let us denote by $p_{\lambda}$ the power-sum symmetric functions. If
$\lambda$ is a partition, let $\lambda_k$ be the number of parts of
size $k$ in $\lambda$ and $f_k(\lambda)$ be the number of fixed points
of a permutation of type $\lambda$. Let $z_{\lambda}$ be the product
over $k$ of $k^{\lambda_k} (\lambda_k)!$, a classical constant
associated to a partition.

\begin{proposition}
  \label{prop_Z_X}
  The symmetric function $Z_{\Lambda  \circ W}$ has the following expression
  \begin{equation}
    \label{Z_CE}
    1+(-t)\sum_{\lambda,|\lambda|\geq 1}(\lambda_1-t)^{\lambda_1-1}\prod_{k \geq
      2}\left( (f_k(\lambda)-t^k)^{\lambda_k}- k \lambda_k (f_k(\lambda)-t^k)^{\lambda_k-1}\right) \frac{p_\lambda}{z_\lambda},
  \end{equation}
  and the symmetric function $Z_X$ has the following expression
   \begin{equation}
     \label{Z_X}
    \sum_{\lambda,|\lambda|\geq 1}(\lambda_1-1)^{\lambda_1-1}\prod_{k \geq
      2}\left( (f_k(\lambda)-1)^{\lambda_k}- k \lambda_k (f_k(\lambda)-1)^{\lambda_k-1}\right) \frac{p_\lambda}{z_\lambda},
  \end{equation}
  where the sums are over the set of non-empty partitions $\lambda$.
\end{proposition}





\bigskip

Recall that it was proved in \cite{anticyclic} the the PreLie operad
is an anticyclic operad. This implies in particular that there is an
action of the symmetric group $\sym_{n+1}$ on the space $\prelie(n)$.
Let us denote by $\widehat{W}$ the corresponding $\sym$-module and let
us compute the symmetric function $Z_{\widehat{W}}$ describing this
action of the symmetric group $\sym_{n+1}$ on the space $\prelie(n)$.

From \cite[Eq. (50)]{anticyclic}, this symmetric function is
characterized by the relation
\begin{equation}
  \label{defi_cyclic}
  1+Z_{\widehat{W}}=p_1(1+Z_{W}+1/Z_{W}).
\end{equation}

\begin{proposition}
  \label{prelie_cyclic}
  The symmetric function $ Z_{\widehat{W}}$ is given by
  \begin{equation}
    \label{Z_Wh}
    \sum_{\lambda,|\lambda|\geq 1,\lambda_1\not=1}(\lambda_1-1)^{\lambda_1-2}\prod_{k \geq
      2}\left( (f_k(\lambda)-1)^{\lambda_k}-k \lambda_k (f_k(\lambda)-1)^{\lambda_k-1}\right) \frac{p_\lambda}{z_\lambda}.
  \end{equation}
\end{proposition}

\begin{proof}
  One has (see \cite{hypertree})
  \begin{equation}
    \left( p_1 \exp\left(-\sum_{k\geq 1} p_k/k\right)\right) \circ Z_W =p_1.
  \end{equation}
  
  Let us introduce new variables:
  \begin{equation}
    y_\ell= p_\ell \circ Z_W,
  \end{equation}
  for $\ell\geq 1$. Then the inverse map is given by
  \begin{equation}
    p_\ell=y_\ell \exp\left(-\sum_k y_{k \ell}/k\right).
  \end{equation}

  Let $\lambda$ be a partition with longest part at most $r$. To
  compute the coefficient of $p_\lambda$ in the symmetric function
  $1+Z_{\widehat{W}}$, it is enough to compute the residue
  \begin{equation}
    \iiint ( 1+Z_{\widehat{W}}) \prod_{i=1}^{r} \frac{dp_i}{p_i^{\lambda_i+1}},
  \end{equation}
  which is equal by formula (\ref{defi_cyclic}) to
  \begin{equation}
    \iiint  p_1(1+Z_W+1/Z_W) \prod_{i=1}^{r} \frac{dp_i}{p_i^{\lambda_i+1}}.
  \end{equation}
  
  One can assume without restriction that all variables $y_j$ and
  $p_j$ vanish if $j>r$.

  We will change the variables to get instead a residue in the
  variables $y$. One has to use the formula
  \begin{equation}
    \prod_{i=1}^{r} dp_i
    =\exp\left(-\sum_{i} \sum_k y_{ik}/k\right)\prod_{i=1}^{r} (1-y_i) dy_i. 
  \end{equation}

  We therefore have to compute the residue
  \begin{equation}
    \iiint y_1 \exp\left(-\sum_{k} y_k/k\right) \left(
      1+y_1+1/y_1 \right) \exp\left(\sum_{i}\lambda_i \sum_k y_{ik}/k\right)\prod_{i=1}^{r} \frac{(1-y_i)}{y_i^{\lambda_i+1}} dy_i. 
  \end{equation}

  Gathering the exponentials and reversing the order of summation, one finds
  \begin{equation}
    \iiint \exp\left(\sum_{k} (f_k(\lambda)-1) y_k/k\right) (1+y_1+y_1^2) \prod_{i=1}^{r} \frac{(1-y_i)}{y_i^{\lambda_i+1}} dy_i. 
  \end{equation}

  This integral decomposes as a product of residues in each 
  variable $y_i$.
  
  Let us discuss first the integral with respect to $y_1$:
  \begin{equation}
    \oint \exp\left( (\lambda_1-1) y_1\right) \frac{(1-y_1^3)}{y_1^{\lambda_1+1}} dy_1. 
  \end{equation}
  This is the sum of two terms:
  \begin{equation}
    \oint \exp\left( (\lambda_1-1) y_1\right)
    \frac{1}{y_1^{\lambda_1+1}} dy_1=
    \frac{(\lambda_1-1)^{\lambda_1}}{\lambda_1 !},
  \end{equation}
  and
  \begin{equation}
    \oint \exp\left( (\lambda_1-1) y_1\right) \frac{(-y_1^3)}{y_1^{\lambda_1+1}} dy_1. =-\frac{(\lambda_1-1)^{\lambda_1-3}}{(\lambda_1-3) !}.
  \end{equation}
  Note that some care must be taken in the second term when
  $\lambda_1 \leq 2$. The resulting sum is zero if
  $\lambda_1=1$ and
  \begin{equation}
     \frac{(\lambda_1-1)^{\lambda_1-2}}{\lambda_1 !}  \text{ if }\lambda_1\not=1.
  \end{equation}

  Let us then discuss the integral with respect to $y_k$ for $k \geq 2$.
  \begin{equation}
    \oint \exp\left( (f_k(\lambda)-1) y_k/k\right) \frac{(1-y_k)}{y_k^{\lambda_k+1}} dy_k.
  \end{equation}
  This is also the sum of two terms 
  \begin{equation}
    (f_k(\lambda)-1)^{\lambda_k}/k^{\lambda_k}/{\lambda_k !}- (f_k(\lambda)-1)^{\lambda_k-1}/k^{\lambda_k-1}/{(\lambda_k-1) !}.
  \end{equation}
  This can be rewritten as
  \begin{equation}
   \frac{(f_k(\lambda)-1)^{\lambda_k}-k \lambda_k
     (f_k(\lambda)-1)^{\lambda_k-1}}{k^{\lambda_k} {\lambda_k !}}.
  \end{equation}

  Gathering all terms and removing the contribution of the empty
  partition gives the result.
\end{proof}

From this, one deduces the following relation.
\begin{theorem}
  One has 
  \begin{equation}
    Z_X-p_1=(p_1 \partial_{p_1} - \operatorname{Id})Z_{\widehat{W}},
  \end{equation}
   \textit{i.e.} the action of $\sym_n$ on $X(n)$ is obtained from
  the action of $\sym_n$ on $\widehat{W}(n)$ by taking the inner
  tensor product with the reflection module of dimension $n-1$ of $\sym_n$.
\end{theorem}

\begin{proof}
  This is a computation done using Prop. \ref{prop_Z_X} and
  \ref{prelie_cyclic}. Indeed, the operator $(p_1 \partial_{p_1} -
  \operatorname{Id})$ acts by multiplication of $p_{\lambda}$ by
  $\lambda_1-1$. There is one subtle point to check, though. In
  formula (\ref{Z_Wh}) for $Z_{\widehat{W}}$, the summation is over
  all partitions of size at least $2$ with $\lambda_1 \not =1$,
  whereas in Formula (\ref{Z_X}) for $Z_X-p_1$, the summation is over
  all partitions of size at least $2$ without further condition. Let
  $\lambda$ be any partition of size at least $2$ with $\lambda_1=1$.
  Let $k$ be the size of the next-to-smallest part of $\lambda$. Then
  $f_k(\lambda)=k\lambda_k+1$ and hence the expression
  \begin{equation}
   (f_k(\lambda)-1)^{\lambda_k}-k \lambda_k
     (f_k(\lambda)-1)^{\lambda_k-1} 
  \end{equation}
  vanishes. It follows that all such partitions do not contribute to
  $Z_X-p_1$ and that one can deduce the expected equation.

  That the operator $(p_1 \partial_{p_1} - \operatorname{Id})$
  corresponds to the inner tensor product by the reflection module is
  a classical fact in the theory of symmetric functions.
\end{proof}

\section{Freeness from homology concentration}

We will use the knowledge of the homology of the bottom line of the
bicomplex to show that the free pre-Lie algebras are free as Lie
algebras.

\subsection{General setting}

\label{general_setting}

Let $P$ be an operad and assume that $P(1)=\kk 1$ and let $P^+$ be the
$\sym$-module such that $P=\kk 1 \oplus P^+$. Let $A$ be a $P$-algebra
in the category of $\sym$-modules.

The structure of $P$-algebra on $A$ is given by a morphism $\mu : P
\circ A \to A$.


Let us define for each $k\geq 0$ a subspace $A_{\geq k}$ of $A$. Let
$A_{\geq 0}$ be $A$. By induction, let $A_{\geq k}$ be the image by
$\mu$ of $P^+ \circ A_{\geq k-1}$.

This is a decreasing filtration of $A$ by subspaces. By construction,
this filtration is in fact a filtration of $P$-algebra.

We will furthermore assume that this filtration is separating, which
is true for instance if $A$ has some auxiliary grading concentrated in
positive degrees.

Let us define $H_0(A)$ to be the degree $0$ component $A_{\geq
  0}/A_{\geq 1}$ of the associated graded $P$-algebra $gr A$.

Let us choose a section of $H_0(A)$ in $A$. Let $\Free_P(H_0(A))$ be
the free $P$-algebra on $H_0(A)$. Then there exists a unique morphism
$\theta$ of $P$-algebra from the free $P$-algebra $\Free_P(H_0(A))$ to
$A$ extending the chosen section.

\begin{proposition}
  The morphism $\theta$ is surjective.
\end{proposition}

\begin{proof}
  Because the filtration is assumed to be separating, it is enough to
  prove that the associated morphism of graded algebras is surjective.

  This is done by induction on the degree associated with the
  filtration. This is true in degree $0$, because the map $\theta$
  comes from a section.

  Let now $[x]$ be a class in $gr A_k$ with $k\geq 1$. Pick a
  representative $x \in A_{\geq k}$ of the class $[x]$. Then by
  hypothesis, $x$ can be written as a sum
  \begin{equation}
    \sum_a \mu(y_a,z_a),
  \end{equation}
  where $y_a \in P^+$ and $z_a\in A_{\geq k-1}$. 



  Then it follows from the fact that we used a filtration of
  $P$-algebra that the class $[x]$ itself can be written as
  \begin{equation}
    \sum_a \mu(y_a,[z_a]),
  \end{equation}
  where $y_a \in P^+$ and $z_a\in A_{\geq k-1}$. Each class $[z_a]$ in
  $gr A_{k-1}$ belongs to the image of $\theta$ by induction.
  Therefore $[x]$ belongs to the image of $\theta$.
\end{proof}

To show that the morphism $\theta$ is an isomorphism, an argument of
equality of dimension (in some appropriate sense) between $A$ and
$\Free_P(H_0(A))$ is therefore sufficient.

\subsection{Application}

Let us apply the general setting above to the case that we are
studying. The operad $P$ is the operad $\lie$. The $P$-algebra $A$ is
the $\sym$-module $W$. As $W$ has no component in degree $0$, one can
apply the previous construction. The space $X$ is exactly the homology
group $H_0(A)$.

\begin{theorem}
  The $\sym$-module $W$ is isomorphic as a Lie algebra in the
  category of $\sym$-modules to the free Lie algebra $\lie \circ X$ on
  the $\sym$-module $X$.
\end{theorem}

\begin{corollary}
  For any vector space $V$, the free pre-Lie algebra on $V$ is
  isomorphic as a Lie algebra to the free Lie algebra on $X(V)$.
\end{corollary}

\begin{proof}  
As follows from section \ref{general_setting}, there is a map $\theta$
from $\lie \circ X$ to $\prelie$ which is surjective. To prove the
theorem, it suffices to compare the dimensions.

Let $f_{X}$ be the generating series 
\begin{equation}
  \label{taylorX}
  \sum_{n\geq  1}(n-1)^{n-1}\frac{x^n}{n!},
\end{equation}
which is associated with the $\sym$-module $X$ of indecomposable
elements of $\prelie$ as a Lie algebra. We have to check that the
generating series of the free Lie algebra $\lie \circ X$ on $X$ is
equal to $f_W$. As the series $f_{\lie}$ is $-\log(1-x)$, this amounts
to the equality
\begin{equation}
  -\log(1-f_X)\stackrel{?}{=}f_W,
\end{equation}
which can be rewritten as 
\begin{equation}
  e^{-f_W}\stackrel{?}{=}1-f_X.
\end{equation}
The constant term in $x$ is $1$ on both sides. Therefore it is enough
to compare the derivatives with respect to $x$:
\begin{equation}
 - f'_W e^{-f_W}\stackrel{?}{=}-f'_X=-(1+x f'_W),
\end{equation}
where the right equality is by comparison of the Taylor expansions
(\ref{taylorW}) and (\ref{taylorX}). But we have by definition $f_W=x
e^{f_W}$, hence
\begin{equation}
  f'_W=e^{f_W}+x f'_W e^{f_W}.
\end{equation}
This proves the expected equality and hence the Theorem.
\end{proof}

\bibliographystyle{alpha}
\bibliography{modulePL}

\end{document}